\documentclass[10pt,reqno]{amsart}
\usepackage{latexsym}
\usepackage{graphicx}
\usepackage{fancyhdr,amssymb}
\usepackage{mathtools}
\DeclarePairedDelimiter{\ceil}{\lceil}{\rceil}

\usepackage{hyperref} 

\title{An Algebraic structure for Square-Prime numbers}
\author{Raghavendra N. Bhat}
\address[Raghavendra N. Bhat]{Department of Mathematics, University of Illinois, 1409 West Green 
Street, Urbana, IL 61801, USA}
\email{rnbhat2@illinois.edu}
\date{March 2023}

\begin{document}

\maketitle

\begin{abstract}
In the paper ``An Abelian Loop for Non-Composites" \cite{Bhat}, we introduced a group-like structure consisting of odd prime numbers and 1, with properties that allowed us to prove analogous results to well known theorems in Number Theory. In this paper, we explore some theorems and conjectures in the SP space.

\end{abstract}
\section{Introduction and Summary of previous results}
\par In our previously published paper titled ``Distribution of Square-Prime Numbers" \cite{BhatSP}, we introduced the concept of SP numbers, which are a product of a prime and a square, with the square not being 1. For example, the first few SP numbers are 8, 12, 18, 20, 27, 28, 32, 44, 45, 48, 50, 52, 63, 68, 72, 75, 76, 80, 92, 98, 99, 108, 112, 116, and 117.
\par The paper mainly proved three results. Firstly, we proved that SP numbers have an asymptotic distribution similar to prime numbers, which means that for a large natural number $n$, the number of SP numbers smaller than $n$ is asymptotic to $\frac{n}{\log n}$. We denoted the number of SP numbers smaller than n by $SP(n)$. The result established was $SP(n) = (\zeta(2) - 1) \cdot \frac{n}{\log n} + O(\frac{n}{\log n^2})$, which is lower than $\pi(n)$ for large $n$.
\par We also proved that there are infinitely many SP numbers with a gap $g$, as long as the gap occurs at least once. As an example, since (27,28) is an SP pair with gap 1, we proved that there are infinitely many $SP$ twins (twins are defined to have gap 1).
\par Thirdly, we examined methods to provide asymptotic approximations for the count of SP numbers that end with a specific digit. For instance, the quantity of SP numbers ending in 1 tends to approach $\frac{1}{400}\frac{n}{\log n}(\zeta(2,1/10) + \zeta(2,9/10) + \zeta(2,3/10) + \zeta(2,7/10) - 4).$
\par In the paper `An Abelian Loop for Non-Composites' \cite{Bhat}, we introduced a group-like structure consisting of 1 and odd prime numbers, allowing us to define a group operation and prove some results from analytic number theory.
\par In this paper, our goal is to study a similar loop for SP numbers, define a group operation and explore the distribution of SP numbers in the group-like structure.
\section{The Loop Q}
Let $Q$ consist of 1 and all SP numbers, with operation $\bullet$, where
for $a, b \in Q, a \bullet b$ gives the smallest element of $Q$ strictly larger than $|a-b|$. This set is an abelian loop because it has the following properties:\begin{enumerate}
\item $Q$ is closed under $\bullet$, because the operation always generates numbers in $Q$, by definition.
\item The $\bullet$ operation is commutative, owing to the absolute value.
Thus, $a \bullet b = b \bullet a$.
\item 1 is the identity element of $Q$ since for all $a \in Q$, $a \bullet 1 = a$.
\item $Q$ is also closed under inversion because for all $a \in S, a \bullet a = 1$. Hence, in $Q$, the inverse of an element is the element itself.
\end{enumerate}
Thus, our infinite set $Q$, consisting of 1 and all SP Numbers is an abelian loop. $Q$ is not a group because there is no associativity. For simplicity, we will define a function $N(x)$ to be the smallest element of $Q$ that is greater than $|x|$.\par Moreover, any $Q_r \subseteq Q$ defined as $Q_r=\{1,\dots,sp_r\}$ (hence containing 1 and all SP numbers upto $sp_r$, the $r^{th}$ SP number) is a sub-loop of $Q$. From \cite{BhatSP} we know that for large $r$, $|Q_r| = O\left(\frac{r}{\log r}\right)$.
\section{Some Results}
We now identify some analogous results from previously proved theorems in SP Numbers \cite{BhatSP}.\\\\
\textbf{Theorem 1:} For all $q \in Q$, there exists $a \in Q$ such that $a \bullet q = a$.\\\\We will first prove the following Lemma.\\
\textbf{Lemma 1:} For any natural number $n$, there exists a sequence of $n$ consecutive natural numbers of which none are SP numbers.
\begin{proof}
Assume for the sake of contradiction that there exists a natural number $k$ such that there is no sequence $n_1, n_2, \dots, n_k$ such that none are SP numbers.\\
$\Rightarrow$ In any given sequence of $k$ consecutive natural numbers, there exists at least one SP number.\\
$\Rightarrow$ For a large natural number $M$, there exist at least $\frac{M}{k}$ SP numbers smaller than $M$, by the pigeon hole principle.\\
Choose $M$ such that $C \log M > k$, where $C$ is a constant, defined later. We thus have $\frac{M}{C\log M}<\frac{M}{k}$. We know from \cite{BhatSP} that $SP(M) = O\left(\frac{M}{\log M}\right)$. Thus, there exists a constant $C$ such that for large $M$, $SP(M) < \frac{M}{C\log M} < \frac{M}{k}$. We thus achieve the required contradiction to our assumption that there exist at least $\frac{M}{k}$ SP numbers smaller than $M$. Thus, for any natural number $n$, there exists a sequence of $n$ consecutive natural numbers of which none are SP numbers.\\
\end{proof}
\noindent\textit{Proof of Theorem 1.}\\
We know from Lemma 1 that there are arbitrarily long natural number strings of non SP numbers. Thus, there is an SP number $a$ that is preceded by at least $2q$ non SPs.
Then $a \bullet q = N(a-q) = a$, because $|a-q|$ is in the string of $2q$ non SPs. $\hfill \square$\\\\
\noindent \textbf{Theorem 2:} For every $n \in \mathbb{N}$, there exist $a_1, a_2, a_3, \dots, a_n \in Q$, all distinct, such that $a_1 \bullet a_2 = a_2 \bullet a_3 = a_3 \bullet a_4 = \dots = a_{n-1} \bullet a_n$.\\\\ We first prove the following Lemma.\\\\
\textbf{Lemma 2:} For any natural number $n$ there exist infinitely many arithmetic progressions of SP numbers having length $n$.\begin{proof}We use the result of Green-Tao’s theorem \cite{greentao} of 2004. We wish to create an arbitrarily long chain of SP Numbers in $Q$ such that each adjacent pair has the same value when the $\bullet$ operation is
performed.\par Green and Tao proved that for any natural number $n$, there exists an arithmetic progression of $n$ primes. Let $p_1, p_2, p_3, \dots, a_n$ be
an arithmetic progression of $n$ primes. Choose an arbitrary square number $r^2$. Multiply each of $p_1, p_2, p_3, \dots, p_n$ by $r^2$ to generate an arithmetic progression of $n$ SP numbers $sp_1 = r^2\cdot p_1, sp_2 = r^2\cdot p_2, \dots, sp_n = r^2\cdot p_n$.\par Since $r^2$ can be arbitrary and there are infinitely many choices for $r^2$, the proof is complete.\\
\end{proof}
\noindent\textit{Proof of Theorem 2.}\\For any $n \in \mathbb{N}$, consider an arithmetic progression $A$ of SP numbers: $sp_1, sp_2, sp_3, ...sp_n$. Let their common difference be $k$. Let $l$ be the smallest SP number larger than $l$. Thus, for every pair $(sp_{j} , sp_{j+1})$ from $A$ we have $sp_{j} \bullet sp_{j+1} = l$. The proof is complete.\\\\
Here is an example of Theorem 2. For $n = 4$, we have the SP arithmetic
progression: $164, 188, 212, 236$. Then, $164 \bullet 188 = 188 \bullet 212 = 212 \bullet 236 = 27.\hfill\square$\\\\
\textbf{Theorem 3:} It is not possible to have $a, b, c \in Q$ such that $a \bullet b = b \bullet c = a \bullet c$, where $a, b, c$ are different from one another.\\\\
\textbf{Lemma 3:} For any $n \in \mathbb{N}$, there exists an SP number between $n$ and $2n$. \begin{proof}Consider arbitrary $n > 8$. We know from Bertrand's postulate \cite{bertrand}, we have a prime $p$ such that $\ceil{\frac{n}{4}}<p<2\ceil{\frac{n}{4}}$. Thus, $n<4p<2n$. Since $4p$ is an SP number, we are done.\\
\end{proof}
\noindent\textit{Proof of Theorem 3.}\\Assume, for the sake of contradiction, that there exists $a,b,c$ such that $a \bullet b = b \bullet c = a \bullet c$.
Without loss of generality, let $a < b < c$. Thus, the smallest SP numbers strictly greater than $|a - b|, |b - c|$ and $|a - c|$ should be the same number. Let this be some SP number $s$. Let $|a - b| = x, |b - c| = y$ and $|a - c| = x + y$.
Thus, we have the nearest SP number to $x, y,$ and $x + y$ equal $s$.
This implies that there is no SP number between $x$ and $x + y$, and between $y$ and $x + y$. But $x + y$ is $\geq 2(min\{x, y\})$. This contradicts Lemma 3 which states that for every number $n \geq 8$, there exists an SP number between $n$ and $2n$.$\hfill\square$\\\\
\textbf{Lemma 4:}. If $t$ is a positive integer, then $N(t)$ and $N(t + 1)$ are either equal or are consecutive in $Q$.
\begin{proof}
Assume, for the sake of contradiction, that $N(t)$ and $N(t+1)$ are neither equal nor consecutive in $Q$. Then, there exists $y \in Q$ such that $y\neq N(t), y\neq N(t+1)$ and $N(t)$ and $y$ are consecutive in $Q$. By definition, $N(t)$ is the smallest SP number greater than $t$. Thus, $N(t)>t$ and for all $k$ between $t$ and $N(t)$, $N(k) = N(t)$, by definition. $\Rightarrow N(t) \geq t+1$. If $N(t)>t+1$, then $N(t+1) = N(t)$. But we assumed that $N(t)$ and $N(t+1)$ are not equal. If $N(t)=t+1$, then $N(t+1)$ is the smallest SP number greater than $N(t)$, i.e. $N(t+1)=y$. This is a contradiction to $y \neq N(t+1)$.
\end{proof}
\noindent \textbf{Theorem 4:} If $a$ and $b$ are SP twins in $S$, then $\forall x \in Q$ such that $x < a$ and $x < b$, either \begin{enumerate}
\item $a \bullet x = b \bullet x$, or
\item $a \bullet x$ and $b \bullet x$ are adjacent in $Q$.
\end{enumerate}
\begin{proof}
If $a$ and $b$ are SP twins, and $x \in Q$, then $b=a+1$, $a \bullet x$ $= N(a-x)$ and $b\bullet x$ $= N(b-x)$. Since $b-x = a - x+ 1$, we are done by Lemma 4.
\end{proof}
\section{Conclusion}
This paper proves some basic results for an algebraic structure consisting of only SP numbers and 1. This gives us some insights on group-like operations in the SP-space.


\begin{thebibliography}{99}


\bibitem{Bhat}
Bhat, R. ``An Abelian Loop for Non-Composite Numbers" Missouri J. Math. Sci. 34 (2). \url{https://arxiv.org/pdf/2110.14716.pdf}

\bibitem{OEIS}
The On-Line Encyclopedia of Integer Sequences, published electronically at
\url{https://oeis.org.} Sequence A228056. \url{https://oeis.org/A228056}

\bibitem{BhatSP}
Bhat, R. ``Distribution of Square-Prime Numbers." Missouri J. Math. Sci. 34 (1) 121 - 126, May 2022. https://doi.org/10.35834/2022/3401121 \url{https://arxiv.org/pdf/2109.10238.pdf}

\bibitem{greentao}
Green, Ben; Tao, Terence (2008).The primes contain arbitrarily long arithmetic
progressions. Annals of Mathematics. 167 (2): 481–547. arXiv:math.NT/0404188

\bibitem{bertrand}
 Dickson, L. E. Bertrand’s Postulate History of the Theory of Numbers, Vol. 1:
Divisibility and Primality. New York: Dover, pp. 435-436, 2005.

\end{thebibliography}
\end{document}